# Stability for the determination of unknown boundary and impedance with a Robin boundary condition

Eva Sincich *

**Abstract**

We consider an inverse problem arising in corrosion detection. We prove a stability result of logarithmic type for the determination of the corroded portion of the boundary and impedance by two measurements on the accessible portion of the boundary.

## 1 Introduction

Let us consider the following boundary value problem

$$\begin{cases} \Delta u = 0 , & \text{in } \Omega , \\ \dfrac{\partial u}{\partial \nu} = g , & \text{on } \Gamma_A , \\ \dfrac{\partial u}{\partial \nu} + \gamma u = 0 , & \text{on } \Gamma_I . \end{cases} \qquad (1.1)$$

where $\Gamma_A$ and $\Gamma_I$ are two open, disjoint portions of $\partial \Omega$ such that $\partial \Omega = \overline{\Gamma_A} \cup \overline{\Gamma_I}$. This problem arises in non-destructive testing and it models the phenomenon of surface corrosion in metals. See [5, 11, 14, 15, 16, 18, 19, 20, 24, 25] for related studies and also [9, 10, 13, 26, 38] for a treatment of the more accurate nonlinear model.

According to this model, $\Omega$ represents the electrostatic conductor, $u$ is the harmonic potential, $g$ is the prescribed current density on the portion of the boundary $\Gamma_A$ accessible to direct inspection. Whereas on $\Gamma_I$, the portion which is out of reach, the potential $u$ satisfies a Robin boundary condition. Such a condition describes the possible presence of corrosion damage on the surface $\Gamma_I$ and the so-called Robin coefficient $\gamma$ is nonnegative and represents the reciprocal of the surface impedance.

The inverse problem that we address here consists in the determination of the unknown Robin coefficient $\gamma$ and the inaccessible part of the boundary $\Gamma_I$ by means of two electrostatic measurements performed on the accessible one $\Gamma_A$.

*Department of Mathematics and Informatics, University of Trieste, Via Valerio, 12/1 34127 Trieste, Italy, esincich@units.it - Department of Mathematics and Informatics, University of Udine, Via delle Scienze 206, 33100 Udine - Italy, eva.sincich@dimi.uniud.it

In the literature, we may find several results concerning the determination of boundaries with homogeneous Dirichlet condition and homogeneous Neumann condition from a single electrostatic measurement (see [4, 6, 7, 12, 30, 31]). On the contrary, as shown in [17, 34] by counterexamples, a single measurement is not sufficient to determine a boundary with a Robin condition.

The global uniqueness issue for the present inverse problem has been solved in [11] for a $C^{2,\alpha}$ boundary. Indeed, the author shows that two Cauchy data pairs, that is $(g, u|_{\Gamma_A}), (\tilde{g}, \tilde{u}|_{\Gamma_A})$ guarantee simultaneously the uniqueness of $\Gamma_I$ and $\gamma$ provided $g$ and $\tilde{g}$ are linearly independent and one of them is positive.

In [32], the authors proved, among various results, the uniqueness issue under the milder regularity assumption of a $C^{1,1}$ domain. Both the above mentioned papers are based on the Martin's integral identity [28] which leads to a contradiction argument. Let us observe that the method used in [11, 32] seems to be not suitable for the quantitative stability issue, which is the aim of the present paper, for this reason we need to introduce a new technique.

In the following, we will prove a stability estimate of logarithmic type for both the Robin coefficient and the surface impedance by two Cauchy data pairs $(g, u|_{\Gamma_A}), (\tilde{g}, \tilde{u}|_{\Gamma_A})$, provided

**i)** $g$ and $\tilde{g}$ are linearly independent and such that, for a given $\kappa > 0$, we have that

$$osc_{\Gamma_A^{2r_0}}\left(\frac{\tilde{g}}{g}\right) = \max_{\Gamma_A^{2r_0}}\left(\frac{\tilde{g}}{g}\right) - \min_{\Gamma_A^{2r_0}}\left(\frac{\tilde{g}}{g}\right) \geqslant \kappa > 0 \tag{1.2}$$

where $\Gamma_A^{2r_0}$ is an inner portion of $\Gamma_A$ which will be described later on.

**ii)** $g$ is positive and such that, for a given $g_0 > 0$, we have that

$$g(x) \geqslant g_0 > 0 \quad \text{for any } x \in \Gamma_A^{2r_0}. \tag{1.3}$$

The stable recovering of the Robin coefficient and the surface impedance needs some a-priori mild assumptions on the two themselves, that is

**i)** the coefficient $\gamma$ is Lipschitz continuous and such that, for a given $\gamma_0 > 0$, we have that

$$\|\gamma\|_{C^{0,1}(\Gamma_I)} \leqslant \gamma_0 \quad , \tag{1.4}$$

**ii)** the domain $\Omega$ is of $C^{1,\alpha}$ class with given bounds, as specified in what follows.

The paper is organized as follows. In Section 2 we introduce the main hypothesis and we formulate our main results Theorem 2.2 and Theorem 2.3. In Section 3, we analyze the direct problem. In Lemma 3.1 we prove a regularity result for the solution $u$ of the direct problem based on the Moser iteration technique. More precisely, we prove that the solution and its first order derivatives are Hölder continuous up to the boundary. In Lemma 3.2 we prove a weak Harnack inequality on the Robin boundary $\Gamma_I$, of the type

$$\rho^{-n}\|u\|_{L^2(\Omega \cap B_{2\rho}(x_0))} \leqslant const. \inf_{x \in \Omega \cap B_\rho(x_0)} u(x) \tag{1.5}$$

where $x_0 \in \Gamma_I$. In Lemma 3.3 we provide a lower bound for the solution $u$ to (1.1) when $g$ is positive and satisfies the condition (1.3). Indeed, we observe that by the Giraud's maximum principle the solution $u$ is positive up to the boundary, namely

$$u(x) \geqslant const. > 0 \text{ in } \overline{\Omega}. \tag{1.6}$$

In order to quantify the positivity of $u$ we combine (1.5) with an iterated use of the interior Harnack inequality obtaining the desired estimate (1.6).

In Section 4 we deal with the inverse problem. We begin by observing that if $u$ and $\tilde{u}$ are two solutions to (1.1) corresponding to current density $g$ and $\tilde{g}$ respectively and $u > 0$, then the function $\lambda = \frac{\tilde{u}}{u}$ is a solution to

$$\begin{cases} \operatorname{div}(u^2 \nabla \lambda) = 0, & \text{in } \Omega, \\ u^2 \dfrac{\partial \lambda}{\partial \nu} = \tilde{g}u - g\tilde{u}, & \text{on } \Gamma_A, \\ u^2 \dfrac{\partial \lambda}{\partial \nu} = 0, & \text{on } \Gamma_I. \end{cases} \tag{1.7}$$

Such a change of the independent variable allows us to treat a new problem with an homogeneous Neumann condition on $\Gamma_I$ which is easier to handle with respect the Robin one. In the following we will denote with $u_1$ and $u_2$ the solution to (1.1) corresponding to domains $\Omega_1$ and $\Omega_2$ (such that $\Gamma_{A,1} = \Gamma_{A,2} = \Gamma_A$), Robin coefficients $\gamma_1$ and $\gamma_2$ and current density $g$ and with $\tilde{u_1}$ and $\tilde{u_2}$ the analogous corresponding to current density $\tilde{g}$. In Proposition 4.2 we provide a smallness control of $\lambda_1 = \frac{\tilde{u_1}}{u_1}$ on the set $\Omega_1 \setminus G$ (see Definition 4.1) by arguments of unique continuation as the three spheres inequality. It means that if the Dirichlet traces of $u_i$ and $\tilde{u}_i$, $i = 1, 2$, are close

$$\|u_1 - u_2\|_{L^2(\Gamma_A)} \leqslant \varepsilon \qquad \|\tilde{u_1} - \tilde{u_2}\|_{L^2(\Gamma_A)} \leqslant \varepsilon \tag{1.8}$$

then

$$\int_{\Omega_1 \setminus G} |\nabla \lambda_1|^2 \leqslant const.(\log |\log(\varepsilon)|)^{-\chi} \tag{1.9}$$

where $\chi > 0$. In Proposition 4.3, due to a further regularity of the boundary of $G$, we give an improvement of the rate of smallness found above. Indeed, the Lipschitz regularity of $\partial G$ allows us to use the cone condition to approach the boundary and to achieve the following estimate

$$\int_{\Omega_1 \setminus G} |\nabla \lambda_1|^2 \leqslant const.(|\log(\varepsilon)|)^{-\eta} \tag{1.10}$$

where $\eta > 0$. In Proposition 4.4 we give a lower bound on the gradient of $\lambda_1$ in the interior of $\Omega_1$, namely

$$\int_{B_\rho(x_0)} |\nabla \lambda_1|^2 \geqslant const. \tag{1.11}$$

where $x_0 \in \Omega_{1,2\rho}$ (see Section 2 for a precise definition). The proof relies in a quantitative evaluation of the following argument. By the linear independence of $g$ and $\tilde{g}$ and by (1.2) we may infer that there exists $z_0 \in \Gamma_A^{2r_0}$ such that

$$\alpha g(z_0) + \beta \tilde{g}(z_0) \geqslant const.\kappa \tag{1.12}$$

for a suitable choice of $\alpha$ and $\beta$. By (1.12) and by unique continuation arguments we observe that $\alpha u_1 + \beta \tilde{u}_1$ is not identically zero in $\Omega_1$, then $\lambda_1$ cannot be constant and hence $|\nabla \lambda_1|$ must be positive in a ball. In Proposition 4.5 we state the following doubling inequality at the boundary

$$\int_{\Omega \cap B_{\beta r}(x_0)} |\nabla \lambda|^2 \leqslant const.\beta^K \int_{\Omega \cap B_r(x_0)} |\nabla \lambda|^2 , \qquad (1.13)$$

with $\beta > 1$ and $x_0 \in \Gamma_I$. Such an inequality combined with the *loglog* smallness control provided in Proposition 4.2 allows us to state a first rough estimate of *loglog* type for $\Omega$ contained in Lemma 4.6, namely

$$d_H(\overline{\Omega_1}, \overline{\Omega_2}) \leqslant const. \log(|\log(\varepsilon)|)^{-\chi} \qquad (1.14)$$

with $\chi > 0$. Consequently in Proposition 4.7 we recall a result obtained in [4], which gives sufficient conditions in order to guarantee that the boundaries of the two $C^{1,\alpha}$ domains $\Omega_1$ and $\Omega_2$ are locally represented as Lipschitz graphs in a common reference system. As a consequence, we notice in Proposition 4.8 that, up to choosing the threshold of the error $\varepsilon$ in (1.8) sufficiently small, the hypothesis of Proposition 4.3 are satisfied. In the proof of Theorem 2.2 we observe that in view of the Lipschitz regularity of the boundary $G$ achieved in Proposition 4.8, the techniques developed in Lemma 4.6 can be carried over by replacing the *loglog* type estimate (1.9) by the *log* type one (1.10), leading to the desired estimate

$$d_H(\overline{\Omega_1}, \overline{\Omega_2}) \leqslant const.|\log(\varepsilon)|^{-\eta} \qquad (1.15)$$

with $\eta > 0$. Finally in the proof of Theorem 2.3 we prove the logarithmic stability for the Robin coefficient $\gamma$ that is

$$\sup_{\substack{P \in \Gamma_{I,1}^{r_0} \\ Q \in B_{2\phi(\varepsilon)}(P) \cap \Gamma_{I,2}^{r_0}}} |\gamma_2(Q) - \gamma_1(P)| \leqslant const.|\log(\varepsilon)|^{-\eta}, \qquad (1.16)$$

with $\eta > 0$. The proof is achieved by combining the stability estimate (1.15) and the arguments developed in [36] relying on quantitative unique continuation techniques.

## 2 The main results

### 2.1 Notations and definitions

We introduce some notations that we shall use in the sequel.
For any $x_0 \in \partial\Omega$ and for any $\rho > 0$ we shall denote

$$\Gamma_A^\rho = \{x \in \Gamma_A : \text{dist}(x, \Gamma_I) > \rho\} \qquad (2.1)$$
$$\Gamma_I^\rho = \{x \in \Gamma_I : \text{dist}(x, \Gamma_A) > \rho\} \qquad (2.2)$$
$$U_A^\rho = \{x \in \overline{\Omega} : \text{dist}(x, \Gamma_A^\rho) < \rho\} \qquad (2.3)$$
$$U_I^\rho = \{x \in \overline{\Omega} : \text{dist}(x, \Gamma_I^\rho) < \rho\} \qquad (2.4)$$
$$\Omega_\rho = \{x \in \Omega : \text{dist}(x, \partial\Omega) > \rho\}. \qquad (2.5)$$

**Definition 2.1.** *Given $\alpha$, $0 < \alpha \leqslant 1$, we shall say that a domain $\Omega$ is of class $C^{1,\alpha}$ with constants $r_0$, $M > 0$ if for any $P \in \Omega$, there exists a rigid transformation of coordinates under which we have $P = 0$ and*

$$\Omega \cap B_{r_0} = \{(x', x_3) : x_3 > \varphi(x')\} \tag{2.6}$$

*where*

$$\varphi : B'_{r_0} \subset \mathbb{R}^{n-1} \to \mathbb{R} \tag{2.7}$$

*is a $C^{1,\alpha}$ function satisfying*

$$|\varphi(0)| = |\nabla\varphi(0)| = 0 \quad \text{and} \quad \|\varphi\|_{C^{1,\alpha}(B'_{r_0})} \leqslant M r_0 , \tag{2.8}$$

*where we denote*

$$\begin{aligned}\|\varphi\|_{C^{1,\alpha}(B'_{r_0})} &= \|\varphi\|_{L^\infty(B'_{r_0})} + r_0 \|\nabla\varphi\|_{L^\infty(B'_{r_0})} + \\ &+ r_0^{1+\alpha} \sup_{\substack{x,y \in B'_{r_0} \\ x \neq y}} \frac{|\nabla\varphi(x) - \nabla\varphi(y)|}{|x-y|^\alpha} .\end{aligned} \tag{2.9}$$

## 2.2 Assumptions and a-priori informations

**Assumption on the domain**

Given $r_0, M > 0$ constants we assume that $\Omega \subset \mathbb{R}^n$ and

$$\Omega \text{ is of } C^{1,\alpha} \text{ class with constants } r_0, M. \tag{2.10}$$

Moreover, we assume that

$$\text{the diameter of } \Omega \text{ is bounded by } d_0 . \tag{2.11}$$

**Assumption on $\gamma$**

Given $\gamma_0 > 0$ constant we assume that the Robin coefficient $\gamma \geqslant 0$ is such that supp $\gamma \subset \Gamma_I$ and

$$\|\gamma\|_{C^{0,1}(\Gamma_I)} \leqslant \gamma_0 . \tag{2.12}$$

**Assumption on $g$ and $\tilde{g}$**

Given $E > 0$ constant we assume that the current fluxes $g$ and $\tilde{g}$ are such that supp $g$, supp $\tilde{g} \subset \Gamma_A$ and

$$\|g\|_{C^{0,\alpha}(\Gamma_A)}, \|\tilde{g}\|_{C^{0,\alpha}(\Gamma_A)} \leqslant E . \tag{2.13}$$

Given $\kappa > 0$ constant we assume that $g$ and $\tilde{g}$ are linearly independent and such that

$$osc_{\Gamma_A^{2r_0}}\left(\frac{\tilde{g}}{g}\right) = \max_{\Gamma_A^{2r_0}}\left(\frac{\tilde{g}}{g}\right) - \min_{\Gamma_A^{2r_0}}\left(\frac{\tilde{g}}{g}\right) \geqslant \kappa > 0. \tag{2.14}$$

Given $g_0 > 0$ constant we assume that $g$ is positive and such that

$$g(x) \geqslant g_0 > 0 \quad \text{for any } x \in \Gamma_A^{2r_0}. \tag{2.15}$$

In the sequel, we shall refer to the a-priori data as the following set of quantities $r_0, M, d_0, \gamma_0, E, \kappa, g_0$.

## 2.3 The main results

**Theorem 2.2.** *Let $\Omega_1, \Omega_2$ be two domains satisfying (2.10) and (2.11). Let $\Gamma_{A,i}, \Gamma_{I,i}$, $i = 1, 2$, be the corresponding accessible and inaccessible parts of the boundaries. Let us assume that $\Gamma_{A,1} = \Gamma_{A,2} = \Gamma_A$. Let $u_i \in H^1(\Omega_i)$ be the solution to (1.1) when $\Omega = \Omega_i$, $\gamma = \gamma_i$ and let $\tilde{u}_i \in H^1(\Omega_i)$ be the solution to (1.1) when $\Omega = \Omega_i$, $\gamma = \gamma_i$, $g = \tilde{g}$. Let (2.12)-(2.15) be satisfied. There exists $\varepsilon_0 > 0$ constant only depending on the a-priori data, such that if for some $\varepsilon, 0 < \varepsilon < \varepsilon_0$ we have*

$$\|u_1 - u_2\|_{L^2(\Gamma_A)} \leqslant \varepsilon \ , \ \|\tilde{u}_1 - \tilde{u}_2\|_{L^2(\Gamma_A)} \leqslant \varepsilon \tag{2.16}$$

*then*

$$d_{\mathcal{H}}(\overline{\Omega_1}, \overline{\Omega_2}) \leqslant \phi(\varepsilon) \ , \tag{2.17}$$

*where $\phi$ is an increasing continuous function on $[0, +\infty)$ which satisfies*

$$\phi(t) \leqslant C|\log(t)|^{-\eta} \tag{2.18}$$

*for every $0 < t < 1$.*

**Theorem 2.3.** *Let the hypothesis of Theorem 2.2 be satisfied. Then, if*

$$\|\tilde{u}_1 - \tilde{u}_2\|_{L^2(\Gamma_A)} \leqslant \varepsilon \ , \ \|\tilde{u}_1 - \tilde{u}_2\|_{L^2(\Gamma_A)} \leqslant \varepsilon \tag{2.19}$$

*we have*

$$\sup_{\substack{P \in \Gamma_{I,1}^{r_0} \\ Q \in B_{2\phi(\varepsilon)}(P) \cap \Gamma_{I,2}^{r_0}}} |\gamma_2(Q) - \gamma_1(P)| \leqslant \phi(\varepsilon) \ , \tag{2.20}$$

*up to a possible replacing of the constants $C$ and $\eta$ in (2.18).*

## 3 The direct problem

**Lemma 3.1** ($C^{1,\alpha}$ regularity up to the boundary). *Let $u \in H^1(\Omega)$ be a solution to (1.1) with $\gamma$ and $g$ satisfying the a-priori assumptions stated above. Then $u \in C^{1,\alpha}(\bar{\Omega})$ and there exists a constant $C > 0$, depending on the a-priori data only, such that*

$$\|u\|_{C^{1,\alpha}(\bar{\Omega})} \leqslant C \ . \tag{3.1}$$

**Proof.** The proof relies on a slight adaptation of the arguments developed in [35, Chap.3] based on the Moser iteration technique [21, Chap.8] and by well-known regularity bounds for the Neumann problem [3, p.667]. □

**Lemma 3.2** (Weak Harnack inequality on the Robin boundary). *Let $u \in H^1(\Omega)$ be a solution to (1.1) with $\gamma$ and $g$ satisfying the a-priori assumptions stated above. Then for every $y_0 \in \Gamma_I^{r_0}$ and for every $0 < \rho < \frac{r_0}{16}$ we have that*

$$\rho^{-n}\|u\|_{L^2(\Omega \cap B_{2\rho}(y_0))} \leqslant C_1 \inf_{y \in \Omega \cap B_\rho(y_0)} u(y) \tag{3.2}$$

*where $C_1$ is a constant depending on the a priori data only.*

**Proof.** Dealing as in [35, Lemma 3.3], we have that

$$\|w\|_{L^{\frac{2\hat{n}}{\hat{n}-2}}(\Omega\cap B_{r_1}(y_0))} \leqslant C'\frac{|\beta+1|+1}{r_2-r_1}\|w\|_{L^2(\Omega\cap B_{r_2}(y_0)))} \quad (3.3)$$

where $0 < r_1 < r_2 \leqslant r_0$, $\beta \in \mathbb{R}\setminus\{0\}$, $\hat{n} = n$ for $n > 2$, $\hat{2} > 2$, $C'$ is a constant depending on the a-priori data only and

$$w = \begin{cases} u^{\frac{\beta+1}{2}} & \text{if } \beta \neq -1 \\ \log(u) & \text{if } \beta = -1 \ . \end{cases}$$

Once that the inequality (3.3) is achieved, we can perform the Moser iteration method arguing as in [21, Theorem 8.18] (see also [35, Lemma 3.3]) in order to obtain the desired weak Harnack inequality (3.2).

□

**Lemma 3.3** (Lower bound on $u$). *Let $u \in H^1(\Omega)$ be a solution to (1.1) with $g \geqslant 0$ satisfying the a priori bounds (2.13) and (2.15). Then, we have that there exists a positive constant $C_0$, depending on the a priori data only, such that*

$$u(x) \geqslant C_0 \text{ for any } x \in \overline{\Omega} \ . \quad (3.4)$$

**Proof.** We observe that there exists a point $x_0 \in \Gamma_I$ such that

$$u(x_0) = \min_{x\in\overline{\Omega}} u(x) = m \ . \quad (3.5)$$

Assume by contradiction that the minimum is not achieved in $\Gamma_I$. Then we would have that, by the maximum principle for harmonic functions [21], there would exist a point $z \in \overline{\Gamma}_A$ such that $m = u(z)$. Being $\frac{\partial u}{\partial \nu} \in C^\alpha(\partial\Omega)$ and being $\frac{\partial u}{\partial \nu} = g \geqslant 0$ on $\Gamma_A$ we have that $\frac{\partial u}{\partial \nu} \geqslant 0$ on $\overline{\Gamma}_A$ which contradicts the Giraud's maximum principle ([22, Theorem 5, pg. 343], see also [27, 29, 23]).
Suppose now that $m \leqslant 0$. Using the Robin condition on $\Gamma_I$ we get that

$$\frac{\partial u}{\partial \nu}(x_0) = -\gamma(x_0)m \geqslant 0 \ , \quad (3.6)$$

which is in contradiction with the Giraud's maximum principle. Hence we deduce that $u(x) > 0$ in $\overline{\Omega}$.
Our purpose now is to obtain a quantitative control of the positivity of $u$ in $\overline{\Omega}$ in terms of the a-priori data only. To this end, we combine the following uniformly boundedness property for harmonic functions (see [21, 35])

$$\sup_{x\in\Omega\cap B_\rho(x_0)} u(x) \leqslant C_2\rho^{-n}\|u\|_{L^2(\Omega\cap B_{2\rho}(x_0))} \ . \quad (3.7)$$

and the (3.2) with $y_0 = x_0$ obtaining

$$\sup_{x\in\Omega\cap B_\rho(z_0)} u(x) \leqslant C_3 u(x_0) \quad (3.8)$$

where $C_2, C_3 > 0$ are constants depending on the a-priori data only.

Let us now choose $y_0 \in \Gamma_A^{2r_0}$ and let $0 < t < \frac{M}{4\sqrt{1+M^2}} r_0$. Using the regularity property $u \in C^{1,\alpha}(\bar{\Omega})$ we have that

$$u(y_0 - t\nu) = u(y_0) + g(z)t + \mathcal{O}(t^{1+\alpha}) \tag{3.9}$$

where $\nu$ is the outward normal to $\Gamma_A$. Recalling that $u(y_0) > 0$ and that $\|u\|_{C^{1,\alpha}(\Gamma_A)} \leqslant C$ we have that choosing $0 < t < \min\{\frac{M}{4\sqrt{1+M^2}} r_0, \left(\frac{g_0}{2C}\right)^{\frac{1}{\alpha}}\}$ the following holds

$$u(y_0 - t\nu) \geqslant \frac{g_0 t}{2} . \tag{3.10}$$

We now consider a point $z_0 \in \Omega \cap B_\rho(x_0)$ such that $B_{\frac{\bar{t}}{8M}}(z_0) \subset \Omega \cap B_\rho(x_0)$. Let us fix $\bar{t} = \min\{\frac{M}{16\sqrt{1+M^2}} r_0, \left(\frac{g_0}{2C}\right)^{\frac{1}{\alpha}}\}$ and let $y_1 = y_0 - \bar{t}\nu$. Let $\gamma$ be a path joining $z_0$ and $y_1$ and let us define $z_i$, $i = 1, \ldots, k$ as follows $z_{i+1} = \gamma(s_i)$, where $t_i = \max\{s : |\gamma(s) - z_i| = \frac{\bar{t}}{4M}\}$ if $|z_0 - z_i| > \frac{\bar{t}}{4M}\}$ otherwise let $i = k$ and stop the process.

We can cover such a path by a chain of finitely many balls $\{B_i\}_{i=1}^N$ with $N \leqslant \frac{16 d_0 nM}{\bar{t}}$ each of which has radius $\frac{\bar{t}}{4M}$ and $B_i \cap B_{i-1} \neq \emptyset$. Then by an iterated use of the Harnack inequality over the chain of balls we obtain that

$$\sup_{x \in B_{\frac{\bar{t}}{4M}}(z_0)} u(x) \geqslant C^N \sup_{x \in B_{\frac{\bar{t}}{4M}}(z_k)} u(x) \tag{3.11}$$

Noticing that $y_1 \in B_{\frac{\bar{t}}{4M}}(z_k)$ we have that (4.9),(3.10) and the above inequality lead

$$u(x_0) \geqslant C g_0 \tag{3.12}$$

where $C$ is constant depending on the a-priori data only. Hence the thesis follows with $C_0 = C g_0$. $\square$

## 4 The inverse problem

We observe that being, by Lemma 3.3, $u_i > 0$, $i = 1, 2$, we can infer that $\lambda_i = \frac{\tilde{u}_i}{u_i}$ is regular in $\Omega_i$. Moreover, after straightforward calculation, we notice that

$$\begin{cases} \text{div}(u_i^2 \nabla \lambda_i) = 0 , & \text{in } \Omega_i , \\ u_i^2 \frac{\partial \lambda_i}{\partial \nu} = \tilde{g} u_i - g \tilde{u}_i , & \text{on } \Gamma_A , \\ u_i^2 \frac{\partial \lambda_i}{\partial \nu} = 0 , & \text{on } \Gamma_{I,i} . \end{cases} \tag{4.1}$$

**Definition 4.1.** *We shall denote with $G$ the connected component of $\Omega_1 \cap \Omega_2$ such that $\Gamma_A \subset \bar{G}$.*

Let us consider $\Omega_1 \setminus \bar{G}$ and let us denote with $N$ the outward unit normal to $\partial(\Omega_1 \setminus \bar{G})$. Then $\lambda_1$ satisfies

$$\begin{cases} \operatorname{div}(u_1^2 \nabla \lambda_1) = 0, & \text{in } \Omega_1 \setminus \bar{G}, \\ u_1^2 \dfrac{\partial \lambda_1}{\partial N} = \eta, & \text{on } \partial(\Omega_1 \setminus \bar{G}) \cap \Gamma_{I,2}, \\ u_1^2 \dfrac{\partial \lambda_1}{\partial N} = 0, & \text{on } \partial(\Omega_1 \setminus \bar{G}) \cap \Gamma_{I,1}, \end{cases} \quad (4.2)$$

where

$$\eta(x) = (u_2)^{-1} \left[ u_1 u_2^2 \left( \frac{\partial \tilde{u}_1}{\partial N} - \frac{\partial \tilde{u}_2}{\partial N} \right) + \frac{\partial \tilde{u}_2}{\partial N} u_1 u_2 (u_2 - u_1) + \tilde{u}_1 \tilde{u}_2^2 \left( \frac{\partial u_1}{\partial N} - \frac{\partial u_2}{\partial N} \right) \right]$$
$$+ (u_2)^{-1} \left[ \frac{\partial u_2}{\partial N} u_1^2 (\tilde{u}_2 - \tilde{u}_1) + \frac{\partial u_2}{\partial N} \tilde{u}_1 (u_1 + u_2)(u_1 - u_2) \right]. \quad (4.3)$$

**Proposition 4.2** (Stability estimate of continuation from Cauchy data). *Let the hypothesis of Theorem 2.2 be satisfied. Then, we have*

$$\int_{\Omega_1 \setminus \bar{G}} |\nabla \lambda_1|^2 \leqslant \omega(\varepsilon), \quad (4.4)$$

*where $\omega$ is an increasing continuous function on $[0, +\infty]$ satisfying*

$$\omega(t) \leqslant C(\log |\log t|)^{-\chi}, \quad (4.5)$$

*where $C, \chi > 0$ are constants depending on the a-priori data only.*

**Proof.** Dealing as in the proof of Proposition 3.1 in [4] (see also [8] and [9, Proposition 3.1]) we can infer that given $P_1 \in \Gamma_A^{2r_0}$, the following estimates hold

$$\|u_1 - u_2\|_{B_{\bar{\rho}}(z_0)}^2 \leqslant C\varepsilon^{2\delta} \quad (4.6)$$

$$\|\nabla u_1 - \nabla u_2\|_{B_{\bar{\rho}}(z_0)}^2 \leqslant C\varepsilon^{2\delta} \quad (4.7)$$

where $C > 0, 0 < \delta < 1$ are constants depending on the a-priori data and $z_0 = P_1 + \frac{M}{4\sqrt{1+M^2}} r_0 \cdot \nu$, $\bar{\rho} = \frac{M r_0}{16\sqrt{1+M^2}}$ with $\nu$ denoting the unit normal at $P_1$. The proof relies on a reformulation of a stability estimate due to Trytten [37] and Payne [33]. Obviously, the same estimates hold true when $u_i$ and $\nabla u_i$, $i = 1, 2$ are replaced by $\tilde{u}_i$ and $\nabla \tilde{u}_i$, $i = 1, 2$.

Following [27], we introduce a regularized distance $\tilde{d}$ from the boundary of $\Omega_1$. We have that there exists $\tilde{d} \in C^2(\Omega_1) \cap C^0(\bar{\Omega}_1)$, satisfying the following properties

**i)** $\gamma_1 \leqslant \dfrac{\operatorname{dist}(x, \partial \Omega_1)}{\tilde{d}(x)} \leqslant \gamma_2$,

**ii)** $|\nabla \tilde{d}(x)| \geqslant c_1$, $\operatorname{dist}(x, \partial \Omega_1) \leqslant b r_0$,

**iii)** $\|\tilde{d}\|_{C^{0,1}} \leqslant c_2 r_0$,

where $\gamma_1, \gamma_2, c_1, c_2, b$ are positive constants depending on $M$ and $\alpha$ only, (see also [4, Lemma 5.2]).

Let us define for every $\rho > 0$

$$\tilde{\Omega}_{1,\rho} = \{x \in \Omega : \tilde{d} > \rho\} . \tag{4.8}$$

It follows that, there exists $a$, $0 < a \leqslant 1$, only depending on $M, \alpha$, such that for every $\rho, 0 < \rho \leqslant ar_0$, $\tilde{\Omega}_{1,\rho}$ is connected with boundary of class $C^1$ and

$$\gamma_1 \rho \leqslant \text{dist}(x, \partial\Omega_1) \leqslant \gamma_2 \rho , \text{ for every } x \in \partial\tilde{\Omega}_{1,\rho} , \tag{4.9}$$

$$|\Omega_1 \setminus \tilde{\Omega}_{1,\rho}| \leqslant \gamma_3 M r_0^{n-1} \rho , \tag{4.10}$$

where $\gamma_3 > 0$ is a constant depending on $M$ and $\alpha$ only. Moreover, for every $x \in \partial\tilde{\Omega}_\rho$, there exists $y \in \partial\Omega$ such that

$$|y - x| = \text{dist}(x, \partial\Omega) , \quad |\nu(x) - \nu(y)| \leqslant \gamma_4 \frac{r^\alpha}{r_0^\alpha} , \tag{4.11}$$

where $\nu(x), \nu(y)$ denote the outer unit normal to $\tilde{\Omega}_1$ at $x$ and to $\Omega$ at $y$ respectively. Moreover, we have also that

$$|\partial\tilde{\Omega}_r|_{n-1} \leqslant \gamma_5 M r_0^{n-1} , \tag{4.12}$$

where $\gamma_5$ is a constant depending on $M$ and $\alpha$ only. Let us define $\theta = \min\{a, \frac{1}{16(1+M^2)\gamma_2}\}$ and let $\bar{r} = r_0\theta$, then we may introduce the set

$$\Gamma_{A,\gamma_2\bar{r}} = \{x \in \Omega_1 \ : \ \text{dist}(x, \partial\Omega_1) = \gamma_2\bar{r}\} . \tag{4.13}$$

We have that

$$\Omega_1 \setminus G \subset [(\Omega_1 \setminus \tilde{\Omega}_{1,\rho}) \setminus G] \cup [\tilde{\Omega}_{1,\rho} \setminus \tilde{V}_r] , \tag{4.14}$$

$$\partial(\tilde{\Omega}_{1,\rho} \setminus \tilde{V}_r) = \tilde{I}_{1,r} \cup \tilde{I}_{2,r} , \tag{4.15}$$

where $\tilde{I}_{1,r}$ is the part of the boundary contained in $\partial\tilde{\Omega}_{1,\rho}$ and $\tilde{I}_{2,r}$ is the part contained in $\partial\tilde{\Omega}_{2,\rho} \cap \partial\tilde{V}_r$. Therefore, we have

$$\int_{\Omega_1 \setminus G} |\nabla\lambda_1|^2 \leqslant \int_{(\Omega_1 \setminus \tilde{\Omega}_{1,\rho}) \setminus G} |\nabla\lambda_1|^2 + \int_{\tilde{\Omega}_{1,\rho} \setminus \tilde{V}_r} |\nabla\lambda_1|^2 . \tag{4.16}$$

By Lemma 3.1 and by Lemma 3.3 we deduce that there exists a positive constant $C$ depending on the a-priori data only such that

$$\|\lambda_1\|_{C^{1,\alpha}(\overline{\Omega}_1)} \leqslant C . \tag{4.17}$$

Hence by (4.17) and (4.10) we have that there exists a constant $C > 0$ depending on the a-priori data only, such that

$$\int_{(\Omega_1 \setminus \tilde{\Omega}_{1,\rho}) \setminus G} |\nabla\lambda_1|^2 \leqslant Cr^\alpha . \tag{4.18}$$

From the divergence theorem we have that

$$\int_{\tilde{\Omega}_{1,\rho}\setminus \tilde{V}_r} |\nabla \lambda_1|^2 \leq (C_0^2)^{-1}\left(\int_{\tilde{I}_{1,r}} |u_1^2 \nabla \lambda_1 \cdot \nu \lambda_1| + \int_{\tilde{I}_{2,r}} |u_1^2 \nabla \lambda_1 \cdot \nu \lambda_1|\right). \quad (4.19)$$

Let $x \in \tilde{I}_{1,r}$. By (4.9), $\text{dist}(x, \partial\Omega) \leq \gamma_2 r$. On the other hand (see [4, Proposition 3.1]) $\text{dist}(x, \Gamma_A) > \gamma_2 r$.
Hence there exists $y \in \partial\Omega_1 \setminus \Gamma_A$ such that $|y - x| = \text{dist}(x, \partial\Omega_1) \leq \gamma_2 r$.
Since $u_1^2 \nabla \lambda_1 \cdot \nu(y) = 0$, by (4.17), (4.9) and by (4.11) we have that there exists a constant $C > 0$ depending on the a-priori data only such that we have

$$|(u_1^2 \nabla \lambda_1 \cdot \nu)(x)| \leq C\left(\frac{r}{r_0}\right)^\alpha. \quad (4.20)$$

Analogously, we have that given $x \in \tilde{I}_{2,r}$ there exists $y \in \partial\Omega_1 \setminus \Gamma_A$ such that $|y - x| = \text{dist}(x, \partial\Omega_2) \leq \gamma_2 r$. Since $(u_2^2 \nabla \lambda_2 \cdot \nu)(y) = 0$, we have that

$$|u_1^2 \nabla \lambda_1 \cdot \nu(x)| \leq \quad u_1^2(x)|\nabla \lambda_1(x) - \nabla \lambda_2(x)| + |u_1^2(x) - u_2^2(x)||\nabla \lambda_2 \cdot \nu(x)| +$$
$$+ u_2^2(x)|\nabla \lambda_2(x) - \nabla \lambda_2(y)|.$$

We notice that

$$\begin{aligned}|\nabla \lambda_1(x) - \nabla \lambda_2(x)| &\leq [(u_1(x)u_2(x))^2]^{-1}\{|(\nabla \tilde{u}_1(x) - \nabla \tilde{u}_2(x))u_1(x)u_2(x)^2| + \\ &\quad |\nabla \tilde{u}_2(x)u_1(x)u_2(x)(u_2(x) - u_1(x))| + \\ &\quad + |(\nabla u_2(x) - \nabla u_1(x))\tilde{u}_1(x)\tilde{u}_2^2(x)| + \\ &\quad + |\nabla u_2(x)u_1^2(x)(\tilde{u}_2(x) - \tilde{u}_1(x))| + \\ &\quad + |\nabla u_2(x)\tilde{u}_1(x)(u_1(x) + u_2(x))(u_1(x) - u_2(x))|\}\end{aligned}$$

Let us define $u = u_2 - u_1$ and $\tilde{u} = \tilde{u}_2 - \tilde{u}_1$. Hence by Lemma 3.1 and Lemma 3.3 we can infer that there exists a positive constant $C$ depending on the a-priori data only such that

$$|u_1^2 \nabla \lambda_1 \cdot \nu(x)| \leq C(r^\alpha + |u(x)| + |\nabla u(x)| + |\tilde{u}(x)| + |\nabla \tilde{u}(x)|). \quad (4.21)$$

Hence by Lemma 3.1, (4.12), (4.16)-(4.21) it follows that

$$\int_{\Omega_1 \setminus G} |\nabla \lambda_1|^2 \leq \quad C(r^\alpha + \max_{\overline{\tilde{V}_r}} |u(x)| + \max_{\overline{\tilde{V}_r}} |\nabla u(x)| +$$
$$+ \max_{\overline{\tilde{V}_r}} |\tilde{u}(x)| + \max_{\overline{\tilde{V}_r}} |\nabla \tilde{u}(x)|) \quad (4.22)$$

By the same arguments discussed in [4, Proposition 3.1] (see also [8] and [9, Theorem 2.1]) and based on an iterated use of the three spheres inequality for solutions to elliptic equations we obtain that there exists $\tau, 0 < \tau < 1, C > 0$ and $r_1 > 0$ such that for any $0 < r < r_1$

$$\int_{B_r(x)} |u|^2 \leq C\left(\int_G |u|^2\right)^{1-\tau^s}\left(\int_{B_r(x)} |u|^2\right)^{\tau^s}, \quad (4.23)$$

$$\int_{B_r(x)} |\nabla u|^2 \leq C\left(\int_G |\nabla u|^2\right)^{1-\tau^s}\left(\int_{B_r(x)} |\nabla u|^2\right)^{\tau^s} \quad (4.24)$$

where $x \in \overline{\tilde{V}_r}$ and $s$ is an integer depending on the a-priori data only. The above estimate are still satisfied when $u$ and $\nabla u$ are replaced by $\tilde{u}$ and $\nabla \tilde{u}$.
Hence we have that by standard estimates for solutions to elliptic equations and by (4.6) and (4.7) we can infer that there exists a positive constant $C$ depending on the a-priori data only such that

$$\int_{B_r(x)} |u|^2 \leqslant C\varepsilon^{2\delta\tau^s} \tag{4.25}$$

$$\int_{B_r(x)} |\nabla u|^2 \leqslant C\varepsilon^{2\delta\tau^s} \tag{4.26}$$

Let us recall the following interpolation inequality

$$\|v\|_{L^\infty(B_\rho)} \leqslant C\left(\left(\int_{B_\rho} v^2\right)^{\frac{\alpha}{2\alpha+n}} \|v\|_{C^{0,\alpha}(B_\rho)}^{\frac{n}{2\alpha+n}} + \frac{1}{\rho^{\frac{n}{2}}}\left(\int_{B_\rho} v^2\right)\right) \tag{4.27}$$

which holds for any function $v$ defined in the ball $B_\rho \subset \mathbb{R}^n$ and for any $\alpha$ such that $0 < \alpha < 1$. By applying (4.27) to $u, \tilde{u}, \nabla u, \nabla \tilde{u}$ in $B_r(x)$ we have that by Lemma 3.1,(4.25) and (4.26) that

$$\|u\|_{L^\infty(B_r)} \leqslant Cr^{-\frac{n}{2}}\varepsilon^{\zeta\tau^s} \tag{4.28}$$

$$\|\nabla u\|_{L^\infty(B_\rho)} \leqslant Cr^{-\frac{n}{2}}\varepsilon^{\zeta\tau^s} \tag{4.29}$$

where $\zeta = \frac{2\alpha\delta}{2\alpha+n}$, $0 < \zeta < 1$ and $C > 0$ depend on the a-priori data only. Inequality (4.28) and (4.29) hold true also for $\tilde{u}, \nabla \tilde{u}$ respectively. Hence by (4.22), (4.28) and (4.29) we obtain that for any $0 < r < r_1$

$$\int_{\Omega_1 \setminus G} |\nabla \lambda_1|^2 \leqslant C\left(r^\alpha + r^{-\frac{n}{2}}\varepsilon^{\zeta\tau^s}\right) \tag{4.30}$$

where $C > 0$ is a constant depending on the a-priori data only. Finally, minimizing the right hand side of (4.30) with respect to $r$ we obtain the thesis. $\square$

**Proposition 4.3** (Improved stability estimate of continuation from Cauchy data). *Let the hypothesis of Proposition 4.2 be fulfilled. In addition, let us assume that there exists a constant $L > 0$ and $r_1, 0 < r_1 < r_0$, such that $\partial G$ is of Lipschitz class with constants $r_1, L$. Then, we have*

$$\int_{\Omega_1 \setminus \bar{G}} |\nabla \lambda_1|^2 \leqslant \phi(\varepsilon) , \tag{4.31}$$

*up to a possible replacing of the constant $C$ and $\eta$ in (2.18).*

**Proof.** We have that by Lemma 3.3 there exists a constant $C > 0$ depending on the a-priori data only such that

$$\int_{\Omega_1 \setminus G} |\nabla \lambda_1|^2 \leqslant C \int_{\partial(\Omega_1 \setminus G)} \lambda_1 (u_1^2 \nabla \lambda_1 \cdot \nu) , \tag{4.32}$$

with $\partial(\Omega_1 \setminus G) \subset (\partial\Omega_1 \setminus \Gamma_A) \cup (\partial\Omega_2 \cap \partial G \setminus U_A^{\frac{r_0}{2}})$.

Since $u_1^2 \nabla \lambda_1 \cdot \nu = 0$ on $\partial\Omega_1 \setminus \Gamma_A$ and $u_2^2 \nabla \lambda_2 \cdot \nu = 0$ on $\partial\Omega_2 \setminus \Gamma_A$ we have by Lemma 3.1 that

$$\int_{\Omega_1 \setminus G} |\nabla \lambda_1|^2 \leqslant C \int_{\partial\Omega_2 \cap \partial G \setminus U_A^{\frac{r_0}{2}}} \lambda_1 [u_1^2 \nabla(\lambda_1 - \lambda_2) \cdot \nu] \leqslant \qquad (4.33)$$

$$\leqslant C_1 \max_{\partial G} |\nabla(\lambda_1 - \lambda_2)| \qquad (4.34)$$

where $C_1 > 0$ is a constant depending on the a-priori data only. By the same argument used in Theorem 4.2 and using the same notations we can deduce that

$$\int_{\Omega_1 \setminus G} |\nabla \lambda_1|^2 \leqslant C_2 (\max_{\partial G} |u| + \max_{\partial G} |\tilde{u}| + \max_{\partial G} |\nabla u| + \max_{\partial G} |\nabla \tilde{u}|) \qquad (4.35)$$

where $C_2 > 0$ is a constant depending on the a-priori data only.

By the Lipschitz regularity of the boundary $\partial G$ it follows that the cone property holds. Precisely, for every point $Q \in \partial G$, there exists a rigid transformation of coordinates under which we have that $Q = 0$ and the finite cone

$$\mathcal{C} = \left\{ x : |x| < r_1, \frac{x \cdot \xi}{|x|} > \cos\theta \right\} \qquad (4.36)$$

with axis in the direction $\xi$ and width $2\theta$, where $\theta = \arctan\frac{1}{L}$, is such that $\mathcal{C} \subset G$.

Let us now consider a point $Q \in \partial G$ and let $Q_0$ be a point lying on the axis $\xi$ of the cone with vertex in $Q = 0$ such that $d_0 = \mathrm{dist}(Q_0, 0) < \frac{r_1}{2}$. Using the notations of Theorem 4.2 we define $\rho_0 = \min\{\bar{\rho}, \frac{r_1}{4}\sin\theta\}$.

By combining the stability estimates near the boundary (4.6) and (4.7) and an iterated use of the three spheres inequality we can claim that

$$\int_{B_{\rho_0}(Q_0)} u^2 \leqslant C \varepsilon^{\delta\tau^s} \qquad (4.37)$$

$$\int_{B_{\rho_0}(Q_0)} |\nabla u|^2 \leqslant C \varepsilon^{\delta\tau^s} \qquad (4.38)$$

where $C > 0$ is a constant depending on the a-priori data only and where the inequalities (4.37) and (4.38) hold also for $\tilde{u}$ and $\nabla \tilde{u}$. By the same argument used in Proposition 3.2 in [4] (see also [9, Theorem 2.1]), based on an iterated use of the three sphere inequality within the cone, we have

$$|u(Q)| + |\tilde{u}(Q)| + |\nabla u(Q)| + |\nabla \tilde{u}(Q)| \leqslant C(|\log \varepsilon|)^{-\eta} \qquad (4.39)$$

where $C > 0, 0 < \eta < 1$ are constants depending on the a priori data only. Finally by (4.35) we obtain the thesis. $\square$

**Proposition 4.4** (Lower bound on the gradient). *Let the hypothesis of Theorem 2.2 be fulfilled. There exist constants $C > 0$ and $\rho_0 > 0$, depending on the a-priori data only, such that, for any $\rho$, $0 < \rho < \rho_0$ and $x_0 \in \Omega_{1,2\rho}$ we have that*

$$\int_{B_\rho(x_0)} |\nabla \lambda_1|^2 \geqslant C\rho^n . \qquad (4.40)$$

**Proof.** Let us consider $x_0 \in \Omega_{1,2\rho}$ and define $\alpha = \tilde{u}_1(x_0)$ and $\beta = -u_1(x_0)$. By the linear independence of $g$ and $\tilde{g}$, we can assume, without loss of generality, that there exists $z_0 \in \Gamma_A^{2r_0}$ such that

$$\alpha g(z_0) + \beta \tilde{g}(z_0) \geqslant c_1 > 0 \tag{4.41}$$

with $c_1 = \dfrac{g_0 \cdot C_0 \cdot \kappa}{2}$ where $C_0$ is the positive constant in (3.4) when $u = u_1$. Infact we may assume that there exists $z_0 \in \Gamma_A^{2r_0}$ such that

$$\alpha g(z_0) + \beta \tilde{g}(z_0) = g(z_0) \cdot \left(\alpha + \beta \frac{\tilde{g}(z_0)}{g(z_0)}\right) \geqslant g(z_0) \cdot \frac{1}{2} osc_{\Gamma_A^{2r_0}} \left(\alpha + \beta \frac{\tilde{g}}{g}\right) =$$

$$= g(z_0) \cdot \frac{1}{2} osc_{\Gamma_A^{2r_0}} \left(\beta \frac{\tilde{g}}{g}\right) = g(z_0) \cdot \frac{1}{2}|\beta| osc_{\Gamma_A^{2r_0}} \left(\frac{\tilde{g}}{g}\right) \geqslant g_0 C_0 \frac{\kappa}{2}.$$

By (2.13) we can infer that for any $z \in B_{\frac{c_1}{2(\alpha+\beta)E}}(z_0) \cap \Gamma_A^{2r_0}$ we have that

$$\alpha g(z) + \beta \tilde{g}(z) \geqslant \frac{c_1}{2} > 0 \,. \tag{4.42}$$

By unique continuation arguments we observe that there exists $y_0 \in B_{\frac{\rho}{2}}(x_0)$ with such that

$$\alpha u_1(y_0) + \beta \tilde{u}_1(y_0) \neq \alpha u_1(x_0) + \beta \tilde{u}_1(x_0) \,. \tag{4.43}$$

By the choice of $\alpha$ and $\beta$ we notice that

$$\alpha u_1(x_0) + \beta \tilde{u}_1(x_0) = 0 \,. \tag{4.44}$$

Being $u_1(x_0), u_1(y_0) > 0$ we have that

$$|\lambda_1(x_0) - \lambda_1(y_0)| = \frac{1}{u_1(y_0)u_1(x_0)} |\alpha u_1(y_0) + \beta \tilde{u}_1(y_0)| \,. \tag{4.45}$$

By the mean value theorem we have that there exists $\xi \in B_{\frac{\rho}{2}}(x_0)$ with $\xi = tx_0 + (1-t)y_0$ for some $0 \leqslant t \leqslant 1$, such that

$$|\nabla \lambda_1(\xi)| = \frac{1}{|x_0 - y_0|} |\lambda_1(x_0) - \lambda_1(y_0)| \,. \tag{4.46}$$

By Lemma 3.1 for $u = u_1$, there exists a positive constant $C$, depending on the a-priori data only, such that

$$\sup_{x \in B_{\frac{\rho}{2}}(x_0)} u_1(x) \leqslant C \tag{4.47}$$

and being $y_0 \in B_{\frac{\rho}{2}}(x_0)$, we have that

$$|\nabla \lambda_1(\xi)| \geqslant \frac{|\alpha u_1(y_0) + \beta \tilde{u}_1(y_0)|}{C^2 \rho} \,. \tag{4.48}$$

By the same argument, based on an iterative use of the Harnack inequality, developed in Lemma 3.3, with $u = \alpha u_1 + \beta \tilde{u}_1$ and with $g$ replaced by $\alpha g + \beta \tilde{g}$, we can infer by (4.42) that there exists a positive constant $c_2$ such that

$$|\nabla \lambda_1(\xi)| \geqslant c_2 \,. \tag{4.49}$$

Being $\lambda_1 = \frac{\tilde{u}_1}{u_1}$ we have that, by (3.1) for $u = \tilde{u}_1$ and $u = u_1$ and by (3.4) for $u = u_1$, there exists a positive constant $C$, depending on the a priori data only, such that

$$\|\lambda_1\|_{C^{1,\alpha}(\bar{\Omega})} \leqslant C. \tag{4.50}$$

Hence we have that

$$|\nabla \lambda_1(x)| \geqslant \frac{c_2}{2} \quad \text{for any} \quad x \in B_{(\frac{c_2}{2C})^{\frac{1}{\alpha}}}(\xi) \ . \tag{4.51}$$

Choosing $\rho_0 = 2(\frac{c_2}{2C})^{\frac{1}{\alpha}}$ and observing that $B_{\frac{\rho}{4}}(\xi) \subset B_\rho(x_0)$ we have that, for any $\rho \leqslant \rho_0$, the following holds

$$\int_{B_\rho(x_0)} |\nabla \lambda_1|^2 \geqslant \int_{B_{\frac{\rho}{4}}(\xi)} |\nabla \lambda_1|^2 \geqslant \rho^n C \tag{4.52}$$

where $C > 0$ is a constant depending on the a priori data only.

$\square$

**Proposition 4.5** (Doubling Inequality at the Boundary). *Let $\lambda \in H^1(\Omega)$ be a solution to*

$$\begin{cases} div(u_1^2 \nabla \lambda) = 0 \ , & in \ \Omega \ , \\ u_1^2 \frac{\partial \lambda}{\partial \nu} = \psi \ , & on \ \Gamma_A \ , \\ u_1^2 \frac{\partial \lambda}{\partial \nu} = 0 \ , & on \ \Gamma_I \end{cases} \tag{4.53}$$

*with $\psi \in L^2(\Gamma_A)$, $\psi \not\equiv 0$ and $\int_{\partial \Omega} \psi = 0$.*
*Let $x_0 \in \Gamma_I$. For every $r > 0$ and every $\beta \geqslant 1$, we have that*

$$\int_{\Omega \cap B_{\beta r}(x_0)} |\nabla \lambda|^2 \leqslant C \beta^K \int_{\Omega \cap B_r(x_0)} |\nabla \lambda|^2 \ , \tag{4.54}$$

*where $C > 0$ and $K > 0$ depend on the a-priori data only.*

**Proof.** For the proof we refer to [2]. The only difference here relies on a more explicit evaluation of the constant $C$ and $K$ in terms of the a-priori data (see [4, 35]). $\square$

**Lemma 4.6** (Loglog stability). *Let the hypothesis of Theorem 2.2 be satisfied. Then, if*

$$\|u_1 - u_2\|_{L^2(\Gamma_A)} \leqslant \varepsilon \ , \|\tilde{u}_1 - \tilde{u}_2\|_{L^2(\Gamma_A)} \leqslant \varepsilon \tag{4.55}$$

*we have*

$$d_{\mathcal{H}}(\partial \Omega_1, \partial \Omega_2) \leqslant \omega(\varepsilon) \ , \tag{4.56}$$

*where $\omega$ is given by (4.5) and $C$ and $a$ are constants depending on the a-priori data only.*

**Proof.** We recall that

$$d = d_{\mathcal{H}}(\partial\Omega_1, \partial\Omega_2) = \max\left(\sup_{x\in\partial\Omega_1} d(x, \partial\Omega_2), \sup_{x\in\partial\Omega_2} d(x, \partial\Omega_1)\right). \tag{4.57}$$

By the Lipschitz regularity of the boundary $\partial\Omega_2$ we have that there exists a point $x_1 \in \Gamma_{I,1}$, such that we have

$$B_{\frac{r_0}{4\sqrt{1+M^2}}}(x_1) \cap \Omega_1 \subset \Omega_1 \setminus \bar{G}. \tag{4.58}$$

Suppose that $r_0 > d$ and set $\beta = \frac{4r_0\sqrt{1+M^2}}{d} > 1$. By Proposition 4.5 we have that

$$\int_{B_{r_0}(x_1)\cap\Omega_1} |\nabla\lambda_1|^2 \leqslant C\beta^K \int_{B_{\frac{d}{4\sqrt{1+M^2}}}(x_1)\cap\Omega_1} |\nabla\lambda_1|^2. \tag{4.59}$$

Since $B_{\frac{d}{4\sqrt{1+M^2}}} \cap \Omega_1 \subset \Omega_1 \setminus \bar{G}$ then by Proposition 4.2, we deduce that

$$\int_{B_{r_0}(x_1)\cap\Omega_1} |\nabla\lambda_1|^2 \leqslant C'\left(\frac{d}{r_0}\right)^{-K} \omega(\varepsilon). \tag{4.60}$$

Combining the lower bound in Proposition 4.4 with an iterated use of the three spheres inequality as in Proposition 4.2 we infer that

$$r_0^n C \leqslant \int_{B_{r_0}(x_1)\cap\Omega_1} |\nabla\lambda_1|^2 \tag{4.61}$$

where $C > 0$ is a constant depending on the a-priori data only.
Hence combining (4.60) and (4.61) we obtain

$$d \leqslant \omega(\varepsilon), \tag{4.62}$$

up to a possible replacing of the constants $C, \chi$ in (4.5).
Suppose now that $r_0 < d$. Of course we have that

$$d \leqslant \frac{d_0}{r_0} r_0. \tag{4.63}$$

Since $B_{\frac{r_0}{4\sqrt{1+M^2}}}(x_1) \cap \Omega_1 \subset \Omega_1 \setminus G$, by the estimate in Proposition 4.2 we have that

$$\int_{B_{\frac{r_0}{4\sqrt{1+M^2}}}(x_1)\cap\Omega_1} |\nabla\lambda_1|^2 \leqslant \omega(\varepsilon). \tag{4.64}$$

Then, by similar arguments as before, we deduce from (4.61) and (4.64) that $r_0 \leqslant \omega(\varepsilon)$. So by (4.63) the conclusion immediately follows. $\square$

**Proposition 4.7** (Graphs condition). *Let $\Omega_1$ and $\Omega_2$ be bounded domains of class $C^{1,\alpha}$ with constants $\rho_0$ and $E$. There exist numbers $\tilde{d}, \tilde{\rho}, \tilde{d} > 0, 0 \leqslant \tilde{\rho}_0 \leqslant \rho_0$ for which the ratios $\frac{\tilde{d}}{\rho_0}$ and $\frac{\tilde{\rho}}{\rho}$ only depend on $\alpha$ and $E$ such that if we have*

$$d_{\mathcal{H}}(\bar{\Omega}_1, \bar{\Omega}_2) \leqslant \tilde{d} \tag{4.65}$$

*then every connected component $G$ of $\Omega_1 \cap \Omega_2$ has a boundary of Lipschitz class with constants $\tilde{\rho}_0, \tilde{E}$ where $\tilde{\rho}_0$ is as above and $\tilde{E} > 0$ only depends on $\alpha$ and $E$.*

**Proof.** For the proof, we refer to [4, Proposition 3.6]. □

**Proposition 4.8.** *Let $u_i, \tilde{u}_i, i = 1, 2$ as in Lemma 4.6. There exists $\varepsilon_0 > 0$ depending on the a-priori data only, such that for any $\varepsilon < \varepsilon_0$ if*

$$\|u_1 - u_2\|_{L^2(\Gamma_A)} \leqslant \varepsilon \ , \|\tilde{u}_1 - \tilde{u}_2\|_{L^2(\Gamma_A)} \leqslant \varepsilon \qquad (4.66)$$

*then every connected component $G$ of $\Omega_1 \cap \Omega_2$ has a boundary of Lipschitz class with constants $\tilde{\rho}_0, L$ depending on the a-priori data only.*

**Proof.** It is well know that in general the Hausdorff distances $d_{\mathcal{H}}(\bar{\Omega}_1, \bar{\Omega}_2)$ and $d_{\mathcal{H}}(\partial \Omega_1, \partial \Omega_2)$ are not equivalent. However, in our regularity assumptions on $\Omega_i, i = 1, 2$, the estimate

$$d_{\mathcal{H}}(\bar{\Omega}_1, \bar{\Omega}_2) \leqslant \omega(\varepsilon) \qquad (4.67)$$

can be derived from (4.56) by arguing as in the proof of Proposition 3.6 [4]. Hence, taking $\varepsilon_0$ small enough so that the right hand side of (4.67) is smaller than $\tilde{d}$ for every $\varepsilon, \varepsilon \leqslant \varepsilon_0$, the result follows from Proposition 4.7.

□

**Proof.** [Theorem 2.2] In view of Proposition 4.8 we have that the hypothesis of Proposition 4.3 are fulfilled. Hence by replacing the rate of stability $\omega$ (defined in Proposition 4.2) with the improved one $\phi$, namely

$$\phi(t) \leqslant C|\log(t)|^{-\eta} \qquad (4.68)$$

in the proof of Proposition 4.6, the thesis follows.

□

**Proof.** [Theorem 2.3] Let us first observe that, in general, the Hausdorff distances $d_H(\overline{\Omega_1}, \overline{\Omega_2})$ and $d_H(\partial \Omega_1, \partial \Omega_2)$ are not equivalent. However, in our regularity assumptions, the following estimate

$$d_H(\partial \Omega_1, \partial \Omega_2) \leqslant \phi(\varepsilon) \qquad (4.69)$$

can be derived from (2.17) using the arguments contained in the proof of Proposition 3.6 in [4]. We consider a point $P \in \Gamma_{I,1}^{r_0}$ and a point $Q \in B_{2\phi(\varepsilon)}(P) \cap \Gamma_{I,2}^{r_0}$. As first step, we notice that being $u_1, u_2 > 0$ in $\overline{\Omega_1}$ and $\overline{\Omega_2}$ respectively, we can compute the difference $\gamma_2(Q) - \gamma_1(P)$ as follows

$$\gamma_2(Q) - \gamma_1(P) = \frac{\partial u_1}{\partial \nu}(P)\frac{1}{u_1(P)} - \frac{\partial u_2}{\partial \nu}(Q)\frac{1}{u_2(Q)} \ . \qquad (4.70)$$

With no loss of generality we may assume that $P, Q \in \overline{\Omega_1}$, hence we have

$$|\gamma_2(Q) - \gamma_1(P)| \leqslant \left|\frac{\partial u_1}{\partial \nu}(P)\frac{1}{u_1(P)} - \frac{\partial u_1}{\partial \nu}(Q)\frac{1}{u_1(Q)}\right| +$$
$$+ \left|\frac{\partial u_2}{\partial \nu}(Q)\frac{1}{u_2(Q)} - \frac{\partial u_1}{\partial \nu}(Q)\frac{1}{u_1(Q)}\right| \qquad (4.71)$$

We can split the first term on the right hand side of (4.71) as follows

$$\left|\frac{\partial u_1}{\partial \nu}(P)\frac{1}{u_1(P)} - \frac{\partial u_1}{\partial \nu}(Q)\frac{1}{u_1(Q)}\right| \leq \left|\frac{\partial u_1}{\partial \nu}(P)\frac{1}{u_1(P)} - \frac{\partial u_1}{\partial \nu}(Q)\frac{1}{u_1(P)}\right| + $$
$$+ \left|\frac{\partial u_1}{\partial \nu}(Q)\frac{1}{u_1(P)} - \frac{\partial u_1}{\partial \nu}(Q)\frac{1}{u_1(Q)}\right|$$

From Lemma 3.1 and Lemma 3.3 we can infer that

$$\left|\frac{\partial u_1}{\partial \nu}(P)\frac{1}{u_1(P)} - \frac{\partial u_1}{\partial \nu}(Q)\frac{1}{u_1(Q)}\right| \leq \frac{C}{C_0}|P-Q|^\alpha + \frac{C}{C_0^2}|P-Q| \,. \quad (4.72)$$

Hence by (4.69) we can infer that

$$\left|\frac{\partial u_1}{\partial \nu}(P)\frac{1}{u_1(P)} - \frac{\partial u_1}{\partial \nu}(Q)\frac{1}{u_1(Q)}\right| \leq \phi(\varepsilon) \,, \quad (4.73)$$

up to a possible replacing of the constants $C$ and $\eta$ in (2.18).
Analogously we can split the second term on the right hand side of (4.71) as follows

$$\left|\frac{\partial u_2}{\partial \nu}(Q)\frac{1}{u_2(Q)} - \frac{\partial u_1}{\partial \nu}(Q)\frac{1}{u_1(Q)}\right| \leq \left|\frac{\partial u_2}{\partial \nu}(Q)\frac{1}{u_2(Q)} - \frac{\partial u_1}{\partial \nu}(Q)\frac{1}{u_2(Q)}\right| + $$
$$+ \left|\frac{\partial u_1}{\partial \nu}(Q)\frac{1}{u_2(Q)} - \frac{\partial u_1}{\partial \nu}(Q)\frac{1}{u_1(Q)}\right|$$

From Lemma 3.1 and Lemma 3.3 we can infer that

$$\left|\frac{\partial u_2}{\partial \nu}(Q)\frac{1}{u_2(Q)} - \frac{\partial u_1}{\partial \nu}(Q)\frac{1}{u_1(Q)}\right| \leq \frac{1}{C_0}\left|\frac{\partial u_2}{\partial \nu}(Q) - \frac{\partial u_1}{\partial \nu}(Q)\right| + \frac{C}{C_0^2}|u_1(Q) - u_2(Q)| \,.$$

Dealing as in Theorem 4.2 in [35] we have that

$$\|u_1 - u_2\|_{C^1(\Gamma_{I,2}^{r_0})} \leq \phi(\varepsilon) \,. \quad (4.74)$$

Hence we have that

$$\left|\frac{\partial u_2}{\partial \nu}(Q)\frac{1}{u_2(Q)} - \frac{\partial u_1}{\partial \nu}(Q)\frac{1}{u_1(Q)}\right| \leq \phi(\varepsilon) \,. \quad (4.75)$$

up to a possible replacing of the constants $C$ and $\eta$ in (2.18). Combining (4.73) and (4.75) we obtain

$$|\gamma_2(Q) - \gamma_1(P)| \leq \phi(\varepsilon) \,. \quad (4.76)$$

Being such an estimate independent from $P$ and $Q$ the thesis follows. $\square$

**Acknowledgments**
The author wishes to thank Prof. G. Alessandrini, Prof. V. Bacchelli and Prof. S. Vessella for their thorough reading as well as for their valuable remarks and comments.